# Calderón's reproducing formulas for the Weinstein $L^2$- multiplier operators

Ahmed Saoudi

**Abstract.** The aim of this work is the study of the Weinstein $L^2$- multiplier operators on $\mathbb{R}_+^{d+1}$ and we give for them Calderón's reproducing formulas and best approximation using the theory of Weinstein transform and reproducing kernels.

## Introduction

The Weinstein operator $\Delta_{W,\alpha}^d$ defined on $\mathbb{R}_+^{d+1} = \mathbb{R}^d \times (0, +\infty)$, by

$$\Delta_{W,\alpha}^d = \sum_{j=1}^{d+1} \frac{\partial^2}{\partial x_j^2} + \frac{2\alpha+1}{x_{d+1}} \frac{\partial}{\partial x_{d+1}} = \Delta_d + L_\alpha, \ \alpha > -1/2,$$

where $\Delta_d$ is the Laplacian operator for the $d$ first variables and $L_\alpha$ is the Bessel operator for the last variable defined on $(0, \infty)$ by

$$L_\alpha u = \frac{\partial^2 u}{\partial x_{d+1}^2} + \frac{2\alpha+1}{x_{d+1}} \frac{\partial u}{\partial x_{d+1}}.$$

The Weinstein operator $\Delta_{W,\alpha}^d$ has several applications in pure and applied mathematics, especially in fluid mechanics [3].

The Weinstein transform generalizing the usual Fourier transform, is given for $\varphi \in L_\alpha^1(\mathbb{R}_+^{d+1})$ and $\lambda \in \mathbb{R}_+^{d+1}$, by

$$\mathcal{F}_W^{\alpha,d}(\varphi)(\lambda) = \int_{\mathbb{R}_+^{d+1}} \varphi(x) \Lambda_\alpha^d(x, \lambda) d\mu_{\alpha,d}(x),$$

where $d\mu_{\alpha,d}(x)$ is the measure on $\mathbb{R}_+^{d+1} = \mathbb{R}^d \times (0, +\infty)$ and $\Lambda_\alpha^d$ is the Weinstein kernel given respectively later by 1.1 and 1.4.

Let $m$ be a function in $L^2_\alpha(\mathbb{R}^{d+1}_+)$ and let $\sigma$ be a positive real number. The Weinstein $L^2$-Multiplier operators is defined for smooth functions $\varphi$ on $\mathbb{R}^{d+1}_+$, as

$$\mathcal{T}_{w,m,\sigma}\varphi(x) := \mathcal{F}^{-1}_{W,\alpha}\left(m_\sigma \mathcal{F}_{W,\alpha}(\varphi)\right)(x), \quad x \in \mathbb{R}^{d+1}_+,$$

where the function $m_\sigma$ is given by

$$m_\sigma(x) = m(\sigma x).$$

These operators are a generalization of the multiplier operators $\mathcal{T}_m$ associated with a bounded function $m$ and given by $\mathcal{T}_m(\varphi) = \mathcal{F}^{-1}(m\mathcal{F}(\varphi))$, where $\mathcal{F}(\varphi)$ denotes the ordinary Fourier transform on $\mathbb{R}^n$. These operators made the interest of several Mathematicians and they were generalized in many settings in [1, 2, 4, 10, 12, 13, 14].

The aim of this paper is to study the multiplier operators $\mathcal{T}_{w,m,\sigma}$. In particular, we give Calderón's reproducing formulas using the theory of Weinstein transform and Weinstein convolution and we use the theory of reproducing kernels to give best approximation of these operators and a Calderón's reproducing formulas of the related extermal function.

This paper is organized as follows. In section 2, we recall some basic harmonic analysis results related with the Weinstein operator $\Delta^d_{W,\alpha}$ and we introduce preliminary facts that will be used later.

In section 3, we study the Weinstein $L^2$-multiplier operators $\mathcal{T}_{w,m,\sigma}$ and we give for them a Plancherel formula and pointwise reproducing formulas. Afterward, we give Calderón's reproducing formulas by using the theory of Weinstein transform.

The last section of this paper is devoted to giving best approximation for the operators $\mathcal{T}_{w,m,\sigma}$ and good estimates of the associated extermal function.

## 1. Harmonic analysis Associated with the Weinstein Operator

In this section, we shall collect some results and definitions from the theory of the harmonic analysis associated with the Weinstein operator $\Delta^d_{W,\alpha}$. Main references are [7, 8, 9].

In the following we denote by

- $\mathbb{R}^{d+1}_+ = \mathbb{R}^d \times (0, \infty)$.
- $x = (x_1, ..., x_d, x_{d+1}) = (x', x_{d+1})$.
- $-x = (-x', x_{d+1})$.
- $C_*(\mathbb{R}^{d+1})$, the space of continuous functions on $\mathbb{R}^{d+1}$, even with respect to the last variable.
- $S_*(\mathbb{R}^{d+1})$, the space of the $C^\infty$ functions, even with respect to the last variable, and rapidly decreasing together with their derivatives.

- $L^p_\alpha(\mathbb{R}^{d+1}_+)$, $1 \leq p \leq \infty$, the space of measurable functions $f$ on $\mathbb{R}^{d+1}_+$ such that

$$\|f\|_{\alpha,p} = \left(\int_{\mathbb{R}^{d+1}_+} |f(x)|^p \, d\mu_\alpha(x)\right)^{1/p} < \infty, \ p \in [1,\infty),$$

$$\|f\|_{\alpha,\infty} = \operatorname{ess\ sup}_{x \in \mathbb{R}^{d+1}_+} |f(x)| < \infty,$$

where

$$d\mu_\alpha(x) = \frac{x_{d+1}^{2\alpha+1}}{(2\pi)^d 2^{2\alpha}\Gamma^2(\alpha+1)} dx. \qquad (1.1)$$

- $\mathcal{A}_\alpha(\mathbb{R}^{d+1}) = \{\varphi \in L^1_\alpha(\mathbb{R}^{d+1}_+); \ \mathcal{F}_{W,\alpha}\varphi \in L^1_\alpha(\mathbb{R}^{d+1}_+)\}$ the Wiener algebra space.

We consider the Weinstein operator $\Delta^d_{W,\alpha}$ defined on $\mathbb{R}^{d+1}_+$ by

$$\Delta^d_{W,\alpha} = \sum_{j=1}^{d+1} \frac{\partial^2}{\partial x_j^2} + \frac{2\alpha+1}{x_{d+1}} \frac{\partial}{\partial x_{d+1}} = \Delta_d + L_\alpha, \ \alpha > -1/2, \qquad (1.2)$$

where $\Delta_d$ is the Laplacian operator for the $d$ first variables and $L_\alpha$ is the Bessel operator for the last variable defined on $(0,\infty)$ by

$$L_\alpha u = \frac{\partial^2 u}{\partial x_{d+1}^2} + \frac{2\alpha+1}{x_{d+1}} \frac{\partial u}{\partial x_{d+1}}.$$

The Weinstein operator $\Delta^d_{W,\alpha}$ have remarkable applications in diffrerent branches of mathematics. For instance, they play a role in Fluid Mechanics [3].

### 1.1. The eigenfunction of the Weinstein operator

For all $\lambda = (\lambda_1, ..., \lambda_{d+1}) \in \mathbb{C}^{d+1}$, the system

$$\frac{\partial^2 u}{\partial x_j^2}(x) = -\lambda_j^2 u(x), \quad \text{if } 1 \leq j \leq d$$

$$L_\alpha u(x) = -\lambda_{d+1}^2 u(x), \qquad (1.3)$$

$$u(0) = 1, \quad \frac{\partial u}{\partial x_{d+1}}(0) = 0, \quad \frac{\partial u}{\partial x_j}(0) = -i\lambda_j, \quad \text{if } 1 \leq j \leq d$$

has a unique solution denoted by $\Lambda^d_\alpha(\lambda, .)$, and given by

$$\Lambda^d_\alpha(\lambda, x) = e^{-i\langle x', \lambda' \rangle} j_\alpha(x_{d+1}\lambda_{d+1}) \qquad (1.4)$$

where $x = (x', x_{d+1})$, $\lambda = (\lambda', \lambda_{d+1})$ and $j_\alpha$ is is the normalized Bessel function of index $\alpha$ defined by

$$j_\alpha(x) = \Gamma(\alpha+1) \sum_{k=0}^{\infty} \frac{(-1)^k x^{2k}}{2^k k! \Gamma(\alpha+k+1)}.$$

The function $(\lambda, x) \mapsto \Lambda^d_\alpha(\lambda, x)$ has a unique extension to $\mathbb{C}^{d+1} \times \mathbb{C}^{d+1}$, and satisfied the following properties.

**Proposition 1.1.** *i). For all $(\lambda, x) \in \mathbb{C}^{d+1} \times \mathbb{C}^{d+1}$ we have*
$$\Lambda_\alpha^d(\lambda, x) = \Lambda_\alpha^d(x, \lambda). \tag{1.5}$$
*ii). For all $(\lambda, x) \in \mathbb{C}^{d+1} \times \mathbb{C}^{d+1}$ we have*
$$\Lambda_\alpha^d(\lambda, -x) = \Lambda_\alpha^d(-\lambda, x). \tag{1.6}$$
*iii). For all $(\lambda, x) \in \mathbb{C}^{d+1} \times \mathbb{C}^{d+1}$ we get*
$$\Lambda_\alpha^d(\lambda, 0) = 1. \tag{1.7}$$
*vi). For all $\nu \in \mathbb{N}^{d+1}$, $x \in \mathbb{R}^{d+1}$ and $\lambda \in \mathbb{C}^{d+1}$ we have*
$$\left| D_\lambda^\nu \Lambda_\alpha^d(\lambda, x) \right| \leq \|x\|^{|\nu|} e^{\|x\| \|\Im \lambda\|} \tag{1.8}$$
*where $D_\lambda^\nu = \partial^\nu / (\partial \lambda_1^{\nu_1} ... \partial \lambda_{d+1}^{\nu_{d+1}})$ and $|\nu| = \nu_1 + ... + \nu_{d+1}$. In particular, for all $(\lambda, x) \in \mathbb{R}^{d+1} \times \mathbb{R}^{d+1}$, we have*
$$\left| \Lambda_\alpha^d(\lambda, x) \right| \leq 1. \tag{1.9}$$

**1.2. The Weinstein transform**

**Definition 1.2.** The Weinstein transform is given for $\varphi \in L_\alpha^1(\mathbb{R}_+^{d+1})$ by
$$\mathcal{F}_{W,\alpha}(\varphi)(\lambda) = \int_{\mathbb{R}_+^{d+1}} \varphi(x) \Lambda_\alpha^d(\lambda, x) d\mu_\alpha(x), \quad \lambda \in \mathbb{R}_+^{d+1}, \tag{1.10}$$
where $\mu_\alpha$ is the measure on $\mathbb{R}_+^{d+1}$ given by the relation (1.1).

Some basic properties of this transform are as follows. For the proofs, we refer [8, 9].

**Proposition 1.3.**
1. *For all $\varphi \in L_\alpha^1(\mathbb{R}_+^{d+1})$, the function $\mathcal{F}_{W,\alpha}(\varphi)$ is continuous on $\mathbb{R}_+^{d+1}$ and we have*
$$\|\mathcal{F}_{W,\alpha} \varphi\|_{\alpha,\infty} \leq \|\varphi\|_{\alpha,1}. \tag{1.11}$$
2. *The Weinstein transform is a topological isomorphism from $\mathcal{S}_*(\mathbb{R}_+^{d+1})$ onto itself. The inverse transform is given by*
$$\mathcal{F}_{W,\alpha}^{-1} \varphi(\lambda) = \mathcal{F}_{W,\alpha} \varphi(-\lambda), \text{ for all } \lambda \in \mathbb{R}_+^{d+1}. \tag{1.12}$$
3. *Parseval formula: For all $\varphi, \phi \in \mathcal{S}_*(\mathbb{R}_+^{d+1})$, we have*
$$\int_{\mathbb{R}_+^{d+1}} \varphi(x) \overline{\phi(x)} d\mu_\alpha(x) = \int_{\mathbb{R}_+^{d+1}} \mathcal{F}_{W,\alpha}(\varphi)(x) \overline{\mathcal{F}_{W,\alpha}(\phi)(x)} d\mu_\alpha(x). \tag{1.13}$$
4. **Plancherel formula**: *For all $\varphi \in \mathcal{S}_*(\mathbb{R}_+^{d+1})$, we have*
$$\|\mathcal{F}_{W,\alpha} \varphi\|_{\alpha,2} = \|\varphi\|_{\alpha,2}. \tag{1.14}$$
5. **Plancherel Theorem**: *The Weinstein transform $\mathcal{F}_{W,\alpha}$ extends uniquely to an isometric isomorphism on $L_\alpha^2(\mathbb{R}_+^{d+1})$.*
6. **Inversion formula**: *Let $\varphi \in L_\alpha^1(\mathbb{R}_+^{d+1})$ such that $\mathcal{F}_{W,\alpha} \varphi \in L_\alpha^1(\mathbb{R}_+^{d+1})$, then we have*
$$\varphi(\lambda) = \int_{\mathbb{R}_+^{d+1}} \mathcal{F}_{W,\alpha} \varphi(x) \Lambda_\alpha^d(-\lambda, x) d\mu_\alpha(x), \text{ a.e. } \lambda \in \mathbb{R}_+^{d+1}. \tag{1.15}$$

## 1.3. The translation operator associated with the Weinstein operator

**Definition 1.4.** The translation operator $\tau_x^\alpha$, $x \in \mathbb{R}_+^{d+1}$ associated with the Weinstein operator $\Delta_{W,\alpha}^d$, is defined for a continuous function $\varphi$ on $\mathbb{R}_+^{d+1}$ which is even with respect to the last variable and for all $y \in \mathbb{R}_+^{d+1}$ by

$$\tau_x^\alpha \varphi(y) = C_\alpha \int_0^\pi \varphi\left(x' + y', \sqrt{x_{d+1}^2 + y_{d+1}^2 + 2x_{d+1}y_{d+1}\cos\theta}\right)(\sin\theta)^{2\alpha}\, d\theta,$$

with

$$C_\alpha = \frac{\Gamma(\alpha+1)}{\sqrt{\pi}\Gamma(\alpha+1/2)}.$$

By using the Weinstein kernel, we can also define a generalized translation, for a function $\varphi \in \mathcal{S}_*(\mathbb{R}^{d+1})$ and $y \in \mathbb{R}_+^{d+1}$ the generalized translation $\tau_x^\alpha \varphi$ is defined by the following relation

$$\mathcal{F}_{W,\alpha}(\tau_x^\alpha \varphi)(y) = \Lambda_\alpha^d(x,y)\mathcal{F}_{W,\alpha}(\varphi)(y). \tag{1.16}$$

The following proposition summarizes some properties of the Weinstein translation operator.

**Proposition 1.5.** *The translation operator $\tau_x^\alpha$, $x \in \mathbb{R}_+^{d+1}$ satisfies the following properties.*

*i). For $\varphi \in \mathbb{C}_*(\mathbb{R}^{d+1})$, we have for all $x, y \in \mathbb{R}_+^{d+1}$*

$$\tau_x^\alpha \varphi(y) = \tau_y^\alpha \varphi(x) \text{ and } \tau_0^\alpha \varphi = \varphi.$$

*ii). Let $\varphi \in L_\alpha^p(\mathbb{R}_+^{d+1})$, $1 \leq p \leq \infty$ and $x \in \mathbb{R}_+^{d+1}$. Then $\tau_x^\alpha \varphi$ belongs to $L_\alpha^p(\mathbb{R}_+^{d+1})$ and we have*

$$\|\tau_x^\alpha \varphi\|_{\alpha,p} \leq \|\varphi\|_{\alpha,p}. \tag{1.17}$$

Note that the $\mathcal{A}_\alpha(\mathbb{R}_+^{d+1})$ is contained in the intersection of $L_\alpha^1(\mathbb{R}_+^{d+1})$ and $L_\alpha^\infty(\mathbb{R}_+^{d+1})$ and hence is a subspace of $L_\alpha^2(\mathbb{R}_+^{d+1})$. For $\varphi \in \mathcal{A}_\alpha(\mathbb{R}_+^{d+1})$ we have

$$\tau_x^\alpha \varphi(y) = C_{\alpha,d} \int_{\mathbb{R}_+^{d+1}} \Lambda_\alpha^d(x,z)\Lambda_\alpha^d(-y,z)\mathcal{F}_{W,\alpha}\varphi(z)d\mu_\alpha(z). \tag{1.18}$$

By using the generalized translation, we define the generalized convolution product $\varphi *_W \psi$ of the functions $\varphi, \psi \in L_\alpha^1(\mathbb{R}_+^{d+1})$ as follows

$$\varphi *_W \psi(x) = \int_{\mathbb{R}_+^{d+1}} \tau_x^\alpha \varphi(-y)\psi(y)d\mu_\alpha(y). \tag{1.19}$$

This convolution is commutative and associative, and it satisfies the following properties.

**Proposition 1.6.** *i) For all $\varphi, \psi \in L_\alpha^1(\mathbb{R}_+^{d+1})$, (resp. $\varphi, \psi \in \mathcal{S}_*(\mathbb{R}_+^{d+1})$), then $\varphi *_W \psi \in L_\alpha^1(\mathbb{R}_+^{d+1})$, (resp. $\varphi *_W \psi \in \mathcal{S}_*(\mathbb{R}_+^{d+1})$) and we have*

$$\mathcal{F}_{W,\alpha}(\varphi *_W \psi) = \mathcal{F}_{W,\alpha}(\varphi)\mathcal{F}_{W,\alpha}(\psi). \tag{1.20}$$

ii) Let $p, q, r \in [1, \infty]$, such that $\frac{1}{p} + \frac{1}{q} - \frac{1}{r} = 1$. Then for all $\varphi \in L^p_\alpha(\mathbb{R}^{d+1}_+)$ and $\psi \in L^q_\alpha(\mathbb{R}^{d+1}_+)$ the function $\varphi *_W \psi$ belongs to $L^r_\alpha(\mathbb{R}^{d+1}_+)$ and we have

$$\|\varphi *_W \psi\|_{\alpha, r} \leq \|\varphi\|_{\alpha, p} \|\psi\|_{\alpha, q}. \tag{1.21}$$

iii) Let $\varphi, \psi \in L^2_\alpha(\mathbb{R}^{d+1}_+)$. Then

$$\varphi *_W \psi = \mathcal{F}^{-1}_{W,\alpha} \left( \mathcal{F}_{W,\alpha}(\varphi) \mathcal{F}_{W,\alpha}(\psi) \right). \tag{1.22}$$

iv) Let $\varphi, \psi \in L^2_\alpha(\mathbb{R}^{d+1}_+)$. Then $\varphi *_W \psi$ belongs to $L^2_\alpha(\mathbb{R}^{d+1}_+)$ if and only if $\mathcal{F}_{W,\alpha}(\varphi)\mathcal{F}_{W,\alpha}(\psi)$ belongs to $L^2_\alpha(\mathbb{R}^{d+1}_+)$ and we have

$$\mathcal{F}_{W,\alpha}(\varphi *_W \psi) = \mathcal{F}_{W,\alpha}(\varphi) \mathcal{F}_{W,\alpha}(\psi). \tag{1.23}$$

v) Let $\varphi, \psi \in L^2_\alpha(\mathbb{R}^{d+1}_+)$. Then

$$\|\varphi *_W \psi)\|_{\alpha, 2} = \|\mathcal{F}_{W,\alpha}(\varphi) \mathcal{F}_{W,\alpha}(\psi)\|_{\alpha, 2}, \tag{1.24}$$

where both sides are finite or infinite.

## 2. Weinstein $L^2$-Multiplier operators on $\mathbb{R}^{d+1}_+$

This section is devoted to study the Weinstein $L^2$-Multiplier operators on $\mathbb{R}^{d+1}_+$ and we establish for them Calderón's reproducing formulas.

**Definition 2.1.** Let $m$ be a function in $L^2_\alpha(\mathbb{R}^{d+1}_+)$ and let $\sigma$ be a positive real number. The Weinstein $L^2$-Multiplier operators is defined for smooth functions $\varphi$ on $\mathbb{R}^{d+1}_+$, as

$$\mathcal{T}_{w,m,\sigma}\varphi(x) := \mathcal{F}^{-1}_{W,\alpha}\left(m_\sigma \mathcal{F}_{W,\alpha}(\varphi)\right)(x), \quad x \in \mathbb{R}^{d+1}_+,$$

where the function $m_\sigma$ is given by

$$m_\sigma(x) = m(\sigma x).$$

According to Proposition 1.6 iii), we can write the operator $\mathcal{T}_{w,m,\sigma}$ as:

$$\mathcal{T}_{w,m,\sigma}\varphi(x) = \mathcal{F}^{-1}_{W,\alpha}(m_\sigma) *_W \varphi(x), \quad x \in \mathbb{R}^{d+1}_+, \tag{2.1}$$

where

$$\mathcal{F}^{-1}_{W,\alpha}(m_\sigma)(x) = \frac{1}{\sigma^{2\alpha+d+2}} \mathcal{F}^{-1}_{W,\alpha}(m)\left(\frac{x}{\sigma}\right).$$

**Proposition 2.2.** (i) For every $m \in L^2_\alpha(\mathbb{R}^{d+1}_+)$, and for every $\varphi \in L^1_\alpha(\mathbb{R}^{d+1}_+)$, the function $\mathcal{T}_{w,m,\sigma}\varphi$ belongs to $L^2_\alpha(\mathbb{R}^{d+1}_+)$, and we have

$$\|\mathcal{T}_{w,m,\sigma}\varphi\|_{\alpha,2} \leq \frac{1}{\sqrt{\sigma^{2\alpha+d+2}}} \|m\|_{\alpha,2} \|\varphi\|_{\alpha,1}.$$

(ii) For every $m \in L^\infty_\alpha(\mathbb{R}^{d+1}_+)$, and for every $\varphi \in L^2_\alpha(\mathbb{R}^{d+1}_+)$, the function $\mathcal{T}_{w,m,\sigma}\varphi$ belongs to $L^2_\alpha(\mathbb{R}^{d+1}_+)$, and we have

$$\|\mathcal{T}_{w,m,\sigma}\varphi\|_{\alpha,2} \leq \|m\|_{\alpha,\infty} \|\varphi\|_{\alpha,2}.$$

(iii) For every $m \in L^2_\alpha(\mathbb{R}^{d+1}_+)$, and for every $\varphi \in L^2_\alpha(\mathbb{R}^{d+1}_+)$, then $\mathcal{T}_{w,m,\sigma}\varphi \in L^\infty_\alpha(\mathbb{R}^{d+1}_+)$, and we have

$$\mathcal{T}_{w,m,\sigma}\varphi(\lambda) = \int_{\mathbb{R}^{d+1}_+} m(\sigma x)\mathcal{F}_{W,\alpha}(\varphi)(x)\Lambda^d_\alpha(-\lambda,x)d\mu_\alpha(x), \quad a.e. \ \lambda \in \mathbb{R}^{d+1}_+$$

and

$$\|\mathcal{T}_{w,m,\sigma}\varphi\|_{\alpha,\infty} \leq \frac{1}{\sqrt{\sigma^{2\alpha+d+2}}}\|m\|_{\alpha,2}\|\varphi\|_{\alpha,2}.$$

*Proof.* i) Let $m \in L^2_\alpha(\mathbb{R}^{d+1}_+)$, and $\varphi \in L^1_\alpha(\mathbb{R}^{d+1}_+)$. From relations (1.11) and (1.14) we get that the function $\mathcal{T}_{w,m,\sigma}\varphi$ belongs to $L^2_\alpha(\mathbb{R}^{d+1}_+)$, and we have

$$\begin{aligned}\|\mathcal{T}_{w,m,\sigma}\varphi\|^2_{\alpha,2} &= \|m_\sigma \mathcal{F}_{W,\alpha}(\varphi)\|^2_{\alpha,2} \\ &\leq \frac{1}{\sqrt{\sigma^{2\alpha+d+2}}}\|m\|^2_{\alpha,2}\|\mathcal{F}_{W,\alpha}(\varphi)\|^2_{\alpha,\infty} \\ &\leq \frac{1}{\sqrt{\sigma^{2\alpha+d+2}}}\|m\|^2_{\alpha,2}\|\varphi\|^2_{\alpha,1}.\end{aligned}$$

ii) The result follows from the Plancherel Theorem for the Weinstein operator.
iii) Let $m \in L^2_\alpha(\mathbb{R}^{d+1}_+)$, and $\varphi \in L^2_\alpha(\mathbb{R}^{d+1}_+)$, then from iversion formula (1.15) we get $\mathcal{T}_{w,m,\sigma}\varphi \in L^\infty_\alpha(\mathbb{R}^{d+1}_+)$, and by relation (1.11) we have

$$\|\mathcal{T}_{w,m,\sigma}\varphi\|_{\alpha,\infty} \leq \|m_\sigma \mathcal{F}_{W,\alpha}(\varphi)\|_{\alpha,1}.$$

Using Hölder's inequality, we get

$$\|\mathcal{T}_{w,m,\sigma}\varphi\|_{\alpha,\infty} \leq \frac{1}{\sqrt{\sigma^{2\alpha+d+2}}}\|m\|_{\alpha,2}\|\varphi\|_{\alpha,2}.$$

□

In the following result, we give Plancherel and pointwise reproducing inversion formulas for the Weinstein $L^2$-Multiplier operators.

**Theorem 2.3.** *Let $m$ be a function in $L^2_\alpha(\mathbb{R}^{d+1}_+)$ satisfying the admissibility condition:*

$$\int_0^\infty |m_\sigma(x)|\frac{d\sigma}{\sigma} = 1, \quad x \in \mathbb{R}^{d+1}_+. \tag{2.2}$$

*i)***Plancherel formula** *For all $\varphi$ in $L^2_\alpha(\mathbb{R}^{d+1}_+)$, we have*

$$\int_{\mathbb{R}^{d+1}_+} |\varphi(x)|^2 d\mu_\alpha(x) = \int_0^\infty \|\mathcal{T}_{w,m,\sigma}\varphi\|^2_{\alpha,2}\frac{d\sigma}{\sigma}.$$

*ii)* **First calderón's formula:** *Let $\varphi$ be a function in $L^1_\alpha(\mathbb{R}^{d+1}_+)$ such that $\mathcal{F}_{W,\alpha}\varphi$ in $L^1_\alpha(\mathbb{R}^{d+1}_+)$ then we have*

$$\varphi(x) = \int_0^\infty \left(\mathcal{T}_{w,m,\sigma}\varphi *_W \mathcal{F}^{-1}_{W,\alpha}(\overline{m_\sigma})\right)(x)\frac{d\sigma}{\sigma}, \quad a.e. \ x \in \mathbb{R}^{d+1}_+.$$

*Proof.* i) From Proposition (1.6) *v*) and relation (2.1) we have

$$\int_0^\infty \int_{\mathbb{R}_+^{d+1}} \|\mathcal{T}_{w,m,\sigma}\varphi\|_{\alpha,2}^2 \frac{d\sigma}{\sigma} = \int_0^\infty \left|\mathcal{F}_{W,\alpha}^{-1}(m_\sigma) *_W \varphi(x)\right|^2 d\mu_\alpha(x) \frac{d\sigma}{\sigma}$$

$$= \int_{\mathbb{R}_+^{d+1}} |\mathcal{F}_{W,\alpha}(x)|^2 \left(\int_0^\infty |m_\sigma|^2 \frac{d\sigma}{\sigma}\right) d\mu_\alpha(x).$$

The result follows from Plancherel Theorem and the assumption (2.2).

ii) Let $\varphi$ be a function in $L_\alpha^1(\mathbb{R}_+^{d+1})$. From Proposition (2.2) i), relation (1.17) and Plancherel Theorem, we have

$$\int_0^\infty \left(\mathcal{T}_{w,m,\sigma}\varphi *_W \mathcal{F}_{W,\alpha}^{-1}(\overline{m_\sigma})\right)(x) \frac{d\sigma}{\sigma} =$$

$$\int_0^\infty \left(\int_{\mathbb{R}_+^{d+1}} \mathcal{T}_{w,m,\sigma}\varphi(y) \overline{\tau_{-x}\left(\mathcal{F}_{W,\alpha}^{-1}(m_\sigma)\right)}(y) d\mu_\alpha(y)\right) \frac{d\sigma}{\sigma} =$$

$$\int_0^\infty \left(\int_{\mathbb{R}_+^{d+1}} \mathcal{F}_{W,\alpha}\left(\mathcal{T}_{w,m,\sigma}\varphi\right)(y) \overline{\mathcal{F}_{W,\alpha}\left(\tau_{-x}\left(\mathcal{F}_{W,\alpha}^{-1}(m_\sigma)\right)\right)}(y) d\mu_\alpha(y)\right) \frac{d\sigma}{\sigma} =$$

$$\int_0^\infty \left(\int_{\mathbb{R}_+^{d+1}} \Lambda_\alpha^d(x,y) \mathcal{F}_{W,\alpha}(\varphi)(y) |m_\sigma(y)|^2 d\mu_\alpha(y)\right) \frac{d\sigma}{\sigma}.$$

By Fubini's theorem, we have

$$\int_0^\infty \left(\mathcal{T}_{w,m,\sigma}\varphi *_W \mathcal{F}_{W,\alpha}^{-1}(\overline{m_\sigma})\right)(x) \frac{d\sigma}{\sigma} =$$

$$\int_{\mathbb{R}_+^{d+1}} \Lambda_\alpha^d(x,y) \mathcal{F}_{W,\alpha}(y) \left(\int_0^\infty |m_\sigma(y)|^2 \frac{d\sigma}{\sigma}\right) d\mu_\alpha(y) =$$

$$\int_{\mathbb{R}_+^{d+1}} \Lambda_\alpha^d(x,y) \mathcal{F}_{W,\alpha}(y) d\mu_\alpha(y) = f(x).$$

□

To establish the Caldersń's reproducing formulas for the Weinstein $L^2$-Multiplier operators, we need the following technical lemma.

**Lemma 2.4.** *Let m be a function in $L_\alpha^2(\mathbb{R}_+^{d+1}) \cap L_\alpha^\infty(\mathbb{R}_+^{d+1})$ satisfy the admissibility condition (2.2). Then the function*

$$\Phi_{\gamma,\delta}(x) = \int_\gamma^\delta |m(\sigma x)|^2 \frac{d\sigma}{\sigma}$$

*belongs to $L_\alpha^2(\mathbb{R}_+^{d+1})$ for all $0 < \gamma < \delta < \infty$ and we have*

$$\Phi_{\gamma,\delta}(x) \in L_\alpha^2(\mathbb{R}_+^{d+1}) \cap L_\alpha^\infty(\mathbb{R}_+^{d+1}).$$

*Proof.* Using Hölder's inequality for the measure $\frac{d\sigma}{\sigma}$, we get

$$\Phi_{\gamma,\delta}(x) \leq \ln(\delta/\gamma) \int_\gamma^\delta |m(\sigma x)|^4 \frac{d\sigma}{\sigma}, \quad x \in \mathbb{R}_+^{d+1}.$$

Wherefore,

$$\begin{aligned}\|\Phi_{\gamma,\delta}\|_{\alpha,2}^2 &\leq \ln(\delta/\gamma)\int_\gamma^\delta\left(\int_{\mathbb{R}_+^{d+1}}|m(\sigma x)|^4 d\mu_\alpha(x)\right)\frac{d\sigma}{\sigma}\\ &\leq \frac{\gamma^{-(2\alpha+d+2)}-\delta^{-(2\alpha+d+2)}}{2\alpha+d+2}\ln(\delta/\gamma)\|m\|_{\alpha,2}^2\|m\|_{\alpha,\infty}^2 < \infty.\end{aligned}$$

On the other hand, from the admissibility condition (2.2), we get

$$\|\Phi_{\gamma,\delta}\|_{\alpha,\infty} \leq 1,$$

which completes the proof. □

**Theorem 2.5.** *(Second Calderón's formula) Let $\varphi \in L_\alpha^2(\mathbb{R}_+^{d+1})$, $m \in L_\alpha^2(\mathbb{R}_+^{d+1}) \cap L_\alpha^\infty(\mathbb{R}_+^{d+1})$ satisfy the admissibility condition (2.2) and $0 < \gamma < \delta < \infty$. Then the function*

$$\varphi_{\gamma,\delta} = \int_\gamma^\delta\left(\mathcal{T}_{w,m,\sigma}\varphi * \mathcal{F}_{W,\alpha}^{-1}(m_\sigma)\right)(x)\frac{d\sigma}{\sigma}, \quad x \in \mathbb{R}_+^{d+1}$$

*belongs to $L_\alpha^2(\mathbb{R}_+^{d+1})$ and satisfies*

$$\lim_{(\gamma,\delta)\to(0,\infty)}\|\varphi_{\gamma,\delta}-\varphi\|_{\alpha,2} = 0. \tag{2.3}$$

*Proof.* According to Prposition (2.2), relation (1.17) and Plancherel Theorem, we obtain

$$\begin{aligned}\varphi_{\gamma,\delta}(x) &= \int_\gamma^\delta\left(\int_{\mathbb{R}_+^{d+1}}\mathcal{T}_{w,m,\sigma}\varphi(y)\overline{\tau_{-x}\left(\mathcal{F}_{W,\alpha}^{-1}(m_\sigma)\right)(y)}d\mu_\alpha(y)\right)\frac{d\sigma}{\sigma}\\ &= \int_\gamma^\delta\left(\int_{\mathbb{R}_+^{d+1}}\Lambda_\alpha^d(x,y)\mathcal{F}_{W,\alpha}(\varphi)(y)|m_\sigma(y)|^2 d\mu_\alpha(y)\right)\frac{d\sigma}{\sigma}.\end{aligned}$$

By FubiniTonnelli's theorem, Hölder's inequality and Lemma (2.4), we get

$$\begin{aligned}&\int_\gamma^\delta\left(\int_{\mathbb{R}_+^{d+1}}\Lambda_\alpha^d(x,y)\mathcal{F}_{W,\alpha}(\varphi)(y)|m_\sigma(y)|^2 d\mu_\alpha(y)\right)\frac{d\sigma}{\sigma}\\ &\leq \int_{\mathbb{R}_+^{d+1}}\mathcal{F}_{W,\alpha}(\varphi)(y)\Phi_{\gamma,\delta}(y)d\mu_\alpha(y)\\ &\leq \|\varphi\|_{\alpha,2}\|\Phi_{\gamma,\delta}\|_{\alpha,2} < \infty.\end{aligned}$$

Then, according to Fubini's theorem and the inversion formula (1.15), we have

$$\begin{aligned}\varphi_{\gamma,\delta}(x) &= \int_{\mathbb{R}_+^{d+1}}\Lambda_\alpha^d(x,y)\mathcal{F}_{W,\alpha}(\varphi)(y)\left(\int_\gamma^\delta|m_\sigma(y)|^2\frac{d\sigma}{\sigma}\right)d\mu_\alpha(y)\\ &= \int_{\mathbb{R}_+^{d+1}}\Lambda_\alpha^d(x,y)\mathcal{F}_{W,\alpha}(\varphi)(y)\Phi_{\gamma,\delta}(y)d\mu_\alpha(y)\\ &= \mathcal{F}_{W,\alpha}^{-1}\left(\mathcal{F}_{W,\alpha}(\varphi)\Phi_{\gamma,\delta}\right)(x).\end{aligned}$$

On the other hand, the function $\Phi_{\gamma,\delta}$ belongs to $L_\alpha^\infty(\mathbb{R}_+^{d+1})$ which allows to see that $\varphi_{\gamma,\delta}$ belongs to $L_\alpha^2(\mathbb{R}_+^{d+1})$ and using Proposition (1.6), we obtain

$$\varphi_{\gamma,\delta} = \mathcal{F}_{W,\alpha}(\varphi)\Phi_{\gamma,\delta}.$$

By the Plancherel formula we get

$$\|\varphi_{\gamma,\delta} - \varphi\|_{\alpha,2}^2 = \int_{\mathbb{R}_+^{d+1}} |\mathcal{F}_{W,\alpha}(\varphi)(y)|^2 (1 - \Phi_{\gamma,\delta}(y))^2 d\mu_\alpha(y).$$

The the admissibility condition (2.2) leads to

$$\lim_{(\gamma,\delta) \to (0,\infty)} \Phi_{\gamma,\delta}(y) = 1, \quad \text{a.e. } y \in \mathbb{R}_+^{d+1}$$

and

$$|\mathcal{F}_{W,\alpha}(\varphi)(y)|^2 (1 - \Phi_{\gamma,\delta}(y))^2 \leq |\mathcal{F}_{W,\alpha}(\varphi)(y)|^2.$$

Finally, the relation (2.3) follows from the dominated convergence theorem. □

In the following section, we study the extremal function associated to the Weinstein $L^2$-multiplier operators.

## 3. The extremal function associated with Weinstein $L^2$-multiplier operators

**Definition 3.1.** Let $\zeta$ be a positive function on $\mathbb{R}_+^{d+1}$ satisfying the following conditions

$$\zeta(\xi) \geq 1, \quad \text{for all} \quad x \in \mathbb{R}_+^{d+1}, \tag{3.1}$$

and

$$\zeta^{-1} \in L_\alpha^1(\mathbb{R}_+^{d+1}). \tag{3.2}$$

We define the Hilbert space $\mathcal{H}_\zeta(\mathbb{R}_+^{d+1})$ by

$$\mathcal{H}_\zeta(\mathbb{R}_+^{d+1}) = \left\{ \varphi \in L_\alpha^2(\mathbb{R}_+^{d+1}) : \sqrt{\zeta}\, \mathcal{F}_{W,\alpha}(\varphi) \in L_\alpha^2(\mathbb{R}_+^{d+1}) \right\},$$

provided with inner product

$$\langle \varphi, \psi \rangle_\zeta = \int_{\mathbb{R}_+^{d+1}} \zeta(\xi) \mathcal{F}_{W,\alpha}(\varphi)(\xi) \overline{\mathcal{F}_{W,\alpha}(\psi)(\xi)} d\mu_\alpha(\xi),$$

and the norm

$$\|\varphi\|_\zeta = \sqrt{\langle \varphi, \varphi \rangle_\zeta}.$$

**Proposition 3.2.** *Let $m$ be a function in $L_\alpha^\infty(\mathbb{R}_+^{d+1})$. Then the Weinstein $L^2$-multiplier operators $\mathcal{T}_{w,m,\sigma}$ are bounded and linear from $\mathcal{H}_\zeta(\mathbb{R}_+^{d+1})$ into $L_\alpha^2(\mathbb{R}_+^{d+1})$ and we have for all $\varphi \in \mathcal{H}_\zeta(\mathbb{R}_+^{d+1})$*

$$\|\mathcal{T}_{w,m,\sigma}\varphi\|_{\alpha,2} \leq \|m\|_{\alpha,\infty} \|\varphi\|_\zeta.$$

*Proof.* Let $\varphi \in \mathcal{H}_\varsigma(\mathbb{R}^{d+1}_+)$. According to Proposition 2.2 (ii), the operator $\mathcal{T}_{w,m,\sigma}$ belongs to $L^2_\alpha(\mathbb{R}^{d+1}_+)$ and we have

$$\|\mathcal{T}_{w,m,\sigma}\varphi\|_{\alpha,2} \leq \|m\|_{\alpha,\infty}\|\varphi\|_{\alpha,2}.$$

On the other hand, by condition (3.1) we have $\|\varphi\|_{\alpha,2} \leq \|\varphi\|_\varsigma$, which gives the result. □

**Definition 3.3.** Let $\eta > 0$ and let $m$ be a function in $L^\infty_\alpha(\mathbb{R}^{d+1}_+)$. We denote by $\langle f, g \rangle_{\varsigma,\eta}$ the inner product defined on the space $\mathcal{H}_\varsigma(\mathbb{R}^{d+1}_+)$ by

$$\langle \varphi, \psi \rangle_{\varsigma,\eta} = \int_{\mathbb{R}^{d+1}_+} \left(\eta\varsigma(\xi) + |m_\sigma(\xi)|^2\right) \mathcal{F}_{W,\alpha}(\varphi)(\xi)\overline{\mathcal{F}_{W,\alpha}(\psi)(\xi)}d\mu_\alpha(\xi), \quad (3.3)$$

and the norm

$$\|\varphi\|_{\varsigma,\eta} = \sqrt{\langle \varphi, \varphi \rangle_{\varsigma,\eta}}.$$

In the following results, we show that the norm $\|\cdot\|_{\varsigma,\eta}$ can be expressed in function of the norm of the Hilbert space $\mathcal{H}_\varsigma(\mathbb{R}^{d+1}_+)$ and the norm of Weinstein $L^2$-Multiplier operators. Moreover, we show the equivalence between the norms $\|\cdot\|_{\varsigma,\eta}$ and $\|\cdot\|_\varsigma$.

**Proposition 3.4.** *Let $m$ be a function in $L^\infty_\alpha(\mathbb{R}^{d+1}_+)$ and $\varphi$ in $\mathcal{H}_\varsigma(\mathbb{R}^{d+1}_+)$*
*(i) The norm $\|\cdot\|_{\varsigma,\eta}$ satisfies:*

$$\|\varphi\|^2_{\varsigma,\eta} = \|\varphi\|^2_\varsigma + \|\mathcal{T}_{w,m,\sigma}\varphi\|^2_{\alpha,2}.$$

*(ii) The norms $\|\cdot\|_{\varsigma,\eta}$ and $\|\cdot\|_\varsigma$ are equivalent and we have*

$$\sqrt{\eta}\,\|\varphi\|_\varsigma \leq \sqrt{\eta + \|m\|^2_{\alpha,\infty}}\,\|\varphi\|_{\varsigma,\eta} \leq \|\varphi\|_\varsigma.$$

*Proof.* (i) follows from the definition of Weinstein $L^2$-Multiplier operators and Plancherel theorem.
(ii) follows from the assertion (i) and Proposition 3.2. □

**Theorem 3.5.** *Let $m$ be a function in $L^\infty_\alpha(\mathbb{R}^{d+1}_+)$. Then the Hilbert space $(\mathcal{H}_\varsigma(\mathbb{R}^{d+1}_+), \langle \cdot, \cdot \rangle_{\varsigma,\eta})$ has the following reproducing Kernel*

$$\Psi_{\varsigma,\eta}(x,y) = \int_{\mathbb{R}^{d+1}_+} \frac{\Lambda^d_\alpha(x,\xi)\Lambda^d_\alpha(-y,\xi)}{\eta\varsigma(\xi) + |m_\sigma(\xi)|^2}d\mu_\alpha(\xi), \quad (3.4)$$

*such that*
*(i) For all $y \in \mathbb{R}^{d+1}_+$, the function $x \mapsto \Psi_{\varsigma,\eta}(x,y)$ belongs to $\mathcal{H}_\varsigma(\mathbb{R}^{d+1}_+)$.*
*(ii) For all $\varphi \in \mathcal{H}_\varsigma(\mathbb{R}^{d+1}_+)$. and $y \in \mathbb{R}^{d+1}_+$, we have the reproducing property*

$$\langle \varphi, \Psi_{\varsigma,\eta}(\cdot,y) \rangle_{\varsigma,\eta} = \varphi(y).$$

*(iii) The Hilbert space $(\mathcal{H}_\varsigma(\mathbb{R}^{d+1}_+), \langle \cdot, \cdot \rangle_\varsigma)$ has the following reproducing Kernel*

$$\Psi_\varsigma(x,y) = \int_{\mathbb{R}^{d+1}_+} \frac{\Lambda^d_\alpha(x,\xi)\Lambda^d_\alpha(-y,\xi)}{\varsigma(\xi)}d\mu_\alpha(\xi). \quad (3.5)$$

*Proof.* (i) Let $y \in \mathbb{R}^{d+1}_+$. From relations (1.9), (3.1) and (3.1), we show that the function
$$f_y : \xi \longrightarrow \frac{\Lambda^d_\alpha(-y,\xi)}{\eta\zeta(\xi) + |m_\sigma(\xi)|^2}$$
belongs to $L^1_\alpha(\mathbb{R}^{d+1}_+) \cap L^2_\alpha(\mathbb{R}^{d+1}_+)$. Hence the function $\Psi_{\zeta,\eta}$ is well defined and by the inversion formula, we botain
$$\Psi_{\zeta,\eta}(x,y) = \mathcal{F}^{-1}_{W,\alpha}(f_y)(x), \quad x \in \mathbb{R}^{d+1}_+.$$

On the other hand, using Plancherel theorem, we get that $\Psi_{\zeta,\eta}(\cdot,y)$ belongs to $L^2_\alpha(\mathbb{R}^{d+1}_+)$ and we have
$$\mathcal{F}_{W,\alpha}\left(\Psi_{\zeta,\eta}(\cdot,y)\right)(\xi) = \frac{\Lambda^d_\alpha(-y,\xi)}{\eta\zeta(\xi) + |m_\sigma(\xi)|^2}, \quad \xi \in \mathbb{R}^{d+1}_+. \tag{3.6}$$

Therefore, by the identity (1.9) we obtain
$$|\mathcal{F}_{W,\alpha}\left(\Psi_{\zeta,\eta}(\cdot,y)\right)(\xi)| \leq \frac{1}{\eta\zeta(\xi)},$$

and
$$\|\Psi_{\zeta,\eta}(\cdot,y)\|^2_\zeta \leq \eta^{-2}\|1/\zeta\|_{\alpha,1} < \infty.$$

This proves that for every $y \in \mathbb{R}^{d+1}_+$, the function $\Psi_{\zeta,\eta}(\cdot,y)$ belongs to $\mathcal{H}_\zeta(\mathbb{R}^{d+1}_+)$.

(ii) Let $\varphi \in \mathcal{H}_\zeta(\mathbb{R}^{d+1}_+)$. and $y \in \mathbb{R}^{d+1}_+$. According to the definition of inner product (3.3) and identity (3.6), we obtain
$$\langle \varphi, \Psi_{\zeta,\eta}(\cdot,y)\rangle_{\zeta,\eta} = \int_{\mathbb{R}^{d+1}_+} \Lambda^d_\alpha(x,\xi)\mathcal{F}_{W,\alpha}(\varphi)(\xi)d\mu_\alpha(\xi).$$

On the other hand, from inequality (3.2) the function $1/\sqrt{\zeta}$ belongs to $L^2_\alpha(\mathbb{R}^{d+1}_+)$. Therefore, the function $\mathcal{F}_{W,\alpha}(\varphi)$ belongs to $L^1_\alpha(\mathbb{R}^{d+1}_+)$ and we have
$$\langle \varphi, \Psi_{\zeta,\eta}(\cdot,y)\rangle_{\zeta,\eta} = \varphi(y).$$

(iii) The result is obtained by taking $m$ a null function and $\eta = 1$. □

The main result of this section can be stated as follows

**Theorem 3.6.** *Let $m$ be a function in $L^\infty_\alpha(\mathbb{R}^{d+1}_+)$ and $\sigma > 0$. For any $h \in L^2_\alpha(\mathbb{R}^{d+1}_+)$ and for any $\eta > 0$, there exists a unique function $\varphi^*_{\eta,h,\sigma}$ where the infimum*
$$\inf_{\varphi \in \mathcal{H}_\zeta}\left\{\eta\|\varphi\|^2_\zeta + \|h - \mathcal{T}_{w,m,\sigma}\varphi\|^2_{\alpha,2}\right\} \tag{3.7}$$
*is attained. Moreover the extremal function $\varphi^*_{\eta,h,\sigma}$ is given by*
$$\varphi^*_{\eta,h,\sigma}(y) = \int_{\mathbb{R}^{d+1}_+} h(x)\overline{\Theta_{\zeta,\eta}(x,y)}d\mu_\alpha(x), \tag{3.8}$$

*where*
$$\Theta_{\zeta,\eta}(x,y) = \int_{\mathbb{R}^{d+1}_+} \frac{m_\sigma(\xi)\Lambda^d_\alpha(x,\xi)}{\eta\zeta(\xi) + |m_\sigma(\xi)|^2}\Lambda^d_\alpha(-y,\xi)d\mu_\alpha(\xi).$$

*Proof.* The existence and unicity of the extremal function $\varphi^*_{\eta,h,\sigma}$ satisfying (3.7) is given by [5, 6, 11]. On the other hand from Theorem 3.5 we have

$$\varphi^*_{\eta,h,\sigma}(y) = \langle h, \mathcal{T}_{w,m,\sigma}(\Psi_{\zeta,\eta}(\cdot,y))\rangle_{\alpha,2}.$$

According to Proposition 2.2 and identity (3.6) we get

$$\begin{aligned}\Theta_{\zeta,\eta}(x,y) &= \mathcal{T}_{w,m,\sigma}(\Psi_{\zeta,\eta}(\cdot,y))(x) \\ &= \int_{\mathbb{R}^{d+1}_+} \frac{m_\sigma(\xi)\Lambda^d_\alpha(x,\xi)}{\eta\zeta(\xi)+|m_\sigma(\xi)|^2}\Lambda^d_\alpha(-y,\xi)d\mu_\alpha(\xi).\end{aligned}$$

□

**Theorem 3.7.** *Let $m$ be a function in $L^\infty_\alpha(\mathbb{R}^{d+1}_+)$ and $h \in L^2_\alpha(\mathbb{R}^{d+1}_+)$. Then the extremal function $\varphi^*_{\eta,h,\sigma}$ satisfies the following properties:*

$$\mathcal{F}_{W,\alpha}(\varphi^*_{\eta,h,\sigma})(\xi) = \frac{\overline{m_\sigma(\xi)}}{\eta\zeta(\xi)+|m_\sigma(\xi)|^2}\mathcal{F}_{W,\alpha}(h)(\xi), \quad \xi \in \mathbb{R}^{d+1}_+,$$

*and*

$$\|\varphi^*_{\eta,h,\sigma}\|^2_\zeta \leq \frac{1}{4\eta}\|h\|^2_{\alpha,2}.$$

*Proof.* Let $y \in \mathbb{R}^{d+1}_+$, then the function

$$g_y : \xi \longrightarrow \frac{m_\sigma(\xi)\Lambda^d_\alpha(-y,\xi)}{\eta\zeta(\xi)+|m_\sigma(\xi)|^2}$$

belongs to $L^1_\alpha(\mathbb{R}^{d+1}_+) \cap L^2_\alpha(\mathbb{R}^{d+1}_+)$ and by the inversion formula we obtain

$$\Theta_{\zeta,\eta}(x,y) = \mathcal{F}^{-1}_{W,\alpha}(g_y)(x), \quad x \in \mathbb{R}^{d+1}_+.$$

Therefore, using Plancherel formula, we have $\Theta_{\zeta,\eta}(\cdot,y)$ belongs to $L^2_\alpha(\mathbb{R}^{d+1}_+)$ and

$$\begin{aligned}\varphi^*_{\eta,h,\sigma}(y) &= \int_{\mathbb{R}^{d+1}_+} \mathcal{F}_{W,\alpha}(h)(\xi)\overline{g_y(\xi)}d\mu_\alpha(\xi) \\ &= \int_{\mathbb{R}^{d+1}_+} \frac{\overline{m_\sigma(\xi)}\mathcal{F}_{W,\alpha}(h)(\xi)}{\eta\zeta(\xi)+|m_\sigma(\xi)|^2}\Lambda^d_\alpha(y,\xi)d\mu_\alpha(\xi).\end{aligned}$$

On the other hand, the function

$$F : \xi \longrightarrow \frac{\overline{m_\sigma(\xi)}\mathcal{F}_{W,\alpha}(h)(\xi)}{\eta\zeta(\xi)+|m_\sigma(\xi)|^2}$$

belongs to $L^1_\alpha(\mathbb{R}^{d+1}_+) \cap L^2_\alpha(\mathbb{R}^{d+1}_+)$ and by the inversion formula we obtain

$$\varphi^*_{\eta,h,\sigma}(y) = \mathcal{F}^{-1}_{W,\alpha}(F)(y).$$

Afterwards, by Plancherel formula, it follows that $\varphi^*_{\eta,h,\sigma}$ belongs to $L^2_\alpha(\mathbb{R}^{d+1}_+)$, and we have

$$\mathcal{F}_{W,\alpha}(\varphi^*_{\eta,h,\sigma})(\xi) = \frac{\overline{m_\sigma(\xi)}}{\eta\zeta(\xi)+|m_\sigma(\xi)|^2}\mathcal{F}_{W,\alpha}(h)(\xi), \quad \xi \in \mathbb{R}^{d+1}_+.$$

Hence

$$|\mathcal{F}_{W,\alpha}(\varphi^*_{\eta,h,\sigma})(\xi)|^2 = \left|\frac{\overline{m_\sigma(\xi)}}{\eta\zeta(\xi)+|m_\sigma(\xi)|^2}\mathcal{F}_{W,\alpha}(h)(\xi)\right|^2,$$

$$\leq \frac{1}{4\eta}|\mathcal{F}_{W,\alpha}(h)(\xi)|^2.$$

Therefore, using Plancherel theorem, we get

$$\|\varphi^*_{\eta,h,\sigma}\|^2_\zeta \leq \frac{1}{4\eta}\|h\|^2_{\alpha,2}.$$

□

**Theorem 3.8.** *(Third Calderón's formula). Let $m$ be a function in $L^\infty_\alpha(\mathbb{R}^{d+1}_+)$ and $\varphi \in \mathcal{H}_\xi(\mathbb{R}^{d+1}_+)$. The extremal function given by*

$$\varphi^*_{\eta,\sigma}(y) = \int_{\mathbb{R}^{d+1}_+} \mathcal{T}_{w,m,\sigma}\varphi(x)\overline{\Theta_{\zeta,\eta}(x,y)}d\mu_\alpha(x). \tag{3.9}$$

*satisfies*

$$\lim_{\eta\to 0^+}\|\varphi^*_{\eta,\sigma}-\varphi\|_\zeta = 0.$$

*Moreover, $\varphi^*_{\eta,\sigma_{\eta>0}}$ converges uniformly to $\varphi$ when $\eta$ converge to $0^+$.*

*Proof.* Let $\varphi \in \mathcal{H}_\xi(\mathbb{R}^{d+1}_+)$, $h = \mathcal{T}_{w,m,\sigma}\varphi$ and $\varphi^*_{\eta,\sigma} = \varphi^*_{\eta,h,\sigma}$. According to Proposition 3.2 the function $h$ belongs to $L^2_\alpha(\mathbb{R}^{d+1}_+)$. From Definition 2.1 and Theorem 3.7, we obtain

$$\mathcal{F}_{W,\alpha}(\varphi^*_{\eta,\sigma})(\xi) = \frac{|m_\sigma(\xi)|^2}{\eta\zeta(\xi)+|m_\sigma(\xi)|^2}\mathcal{F}_{W,\alpha}(\varphi)(\xi), \quad \xi \in \mathbb{R}^{d+1}_+.$$

Hence, it follows that

$$\mathcal{F}_{W,\alpha}(\varphi^*_{\eta,\sigma}-\varphi)(\xi) = \frac{-\eta\zeta(\xi)}{\eta\zeta(\xi)+|m_\sigma(\xi)|^2}\mathcal{F}_{W,\alpha}(\varphi)(\xi), \quad \xi \in \mathbb{R}^{d+1}_+. \tag{3.10}$$

Therefore,

$$\|\varphi^*_{\eta,\sigma}-\varphi\|^2_\zeta = \int_{\mathbb{R}^{d+1}_+}\frac{\eta^2\zeta^3(\xi)|\mathcal{F}_{W,\alpha}(\varphi)(\xi)|^2}{(\eta\zeta(\xi)+|m_\sigma(\xi)|^2)^2}d\mu_\alpha(x).$$

Then, from the dominated convergence theorem and the following inequality

$$\frac{\eta^2\zeta^3(\xi)|\mathcal{F}_{W,\alpha}(\varphi)(\xi)|^2}{(\eta\zeta(\xi)+|m_\sigma(\xi)|^2)^2} \leq \zeta(\xi)|\mathcal{F}_{W,\alpha}(\varphi)(\xi)|^2,$$

we deduce that

$$\lim_{\eta\to 0^+}\|\varphi^*_{\eta,\sigma}-\varphi\|_\zeta = 0.$$

On the other hand, from relation (3.2) the function $1/\sqrt{\zeta}$ belongs to $L^2_\alpha(\mathbb{R}^{d+1}_+)$, hence $\mathcal{F}_{W,\alpha}(\varphi)$ belongs to $L^1_\alpha(\mathbb{R}^{d+1}_+) \cap L^2_\alpha(\mathbb{R}^{d+1}_+)$ for all $\varphi \in \mathcal{H}_\xi(\mathbb{R}^{d+1}_+)$. Then, according to (3.10) and the inversion transform for the Weinstein formula, we get

$$\varphi_{\eta,\sigma}^*(y) - \varphi(y) = \int_{\mathbb{R}_+^{d+1}} \frac{-\eta\zeta(\xi)\mathcal{F}_{W,\alpha}(\varphi)(\xi)}{\eta\zeta(\xi) + |m_\sigma(\xi)|^2} \Lambda(-y,\xi) d\mu_\alpha(x).$$

By using the dominated convergence theorem and the fact

$$\frac{\eta\zeta(\xi)|\mathcal{F}_{W,\alpha}(\varphi)(\xi)|^2}{\eta\zeta(\xi) + |m_\sigma(\xi)|^2} \leq |\mathcal{F}_{W,\alpha}(\varphi)(\xi)|,$$

we deduce that

$$\lim_{\eta \to 0^+} \sup_{y \in \mathbb{R}_+^{d+1}} \|\varphi_{\eta,\sigma}^*(y) - \varphi(y)\| = 0.$$

which completes the proof of the Theorem. □

Ahmed Saoudi
Université de Tunis El Manar, Faculté des sciences de Tunis,
LR13ES06 Laboratoire de Fonctions spéciales, analyse harmonique et analogues, 2092 Tunis, Tunisia
Department of Mathematics, Faculty of Science, Northern Borders University, Arar, Saudi Arabia.
e-mail: `ahmed.saoudi@ipeim.rnu.tn`